\documentclass[fleqn,preprint,3p,a4paper]{elsarticle}

\usepackage{amssymb,mathrsfs,latexsym,amsfonts}
\usepackage{amsmath}
\usepackage{amsthm}
\usepackage{dcolumn}
\usepackage{endnotes}
\usepackage{tabularx}
\usepackage[matrix,arrow]{xy}
\usepackage{wasysym}

\newtheorem{theorem}{Theorem}[section]
\newtheorem{proposition}[theorem]{Proposition}
\newtheorem{lemma}[theorem]{Lemma}
\newtheorem{corollary}[theorem]{Corollary}

\theoremstyle{definition}
\newtheorem{definition}[theorem]{Definition}
\newtheorem{example}[theorem]{Example}
\newtheorem{remark}[theorem]{Remark}
\newtheorem{question}[theorem]{Question}

\journal{Journal}

\begin{document}

\begin{frontmatter}



\title{A supplement on feathered gyrogroups \tnoteref{t1}}
\tnotetext[t1]{This research was supported by the National Natural Science Foundation of China (Nos. 12071199, 11661057), the Natural Science Foundation of Jiangxi Province, China (No. 20192ACBL20045).}


\author[M. Bao]{Meng Bao}
\ead{mengbao95213@163.com}
\address[M. Bao]{College of Mathematics, Sichuan University, Chengdu 610064, China}

\author[X. Ling]{Xuewei Ling}
\ead{781736783@qq.com}
\address[X. Ling]{Institute of Mathematics, Nanjing Normal University, Nanjing 210046, China}

\author[X. Xu]{Xiaoquan Xu\corref{mycorrespondingauthor}}
\cortext[mycorrespondingauthor]{Corresponding author.}
\ead{xiqxu2002@163.com}
\address[X. Xu]{Fujian Key Laboratory of Granular Computing and Applications, Minnan Normal University, Zhangzhou 363000, China}

\begin{abstract}
A topological gyrogroup is a gyrogroup endowed with a topology such that the binary operation is jointly continuous and the inverse mapping is also continuous. It is shown that each compact subset of a topological gyrogroup with an $\omega^{\omega}$-base is metrizable, which deduces that if $G$ is a topological gyrogroup with an $\omega^{\omega}$-base and is a $k$-space, then it is sequential. Moreover, for a feathered strongly topological gyrogroup $G$, based on the characterization of feathered strongly topological gyrogroups, we show that if $G$ has countable $cs^{*}$-character, then it is metrizable; and it is also shown that $G$ has a compact resolution swallowing the compact sets if and only if $G$ contains a compact $L$-subgyrogroup $H$ such that the quotient space $G/H$ is a Polish space.
\end{abstract}

\begin{keyword}
Topological gyrogroups, metrizable, $\omega^{\omega}$-base, feathered
\MSC 22A22; 54A20; 20N05; 18A32.

\end{keyword}




\end{frontmatter}


\section{Introduction}

In the field of Topological Algebra, topological groups are standard researching objects and have been studied for many years, see \cite{AA}. Moreover, the combination of topology and non-associative algebra has attracted the attention of many scholars. For example, in \cite{CZ}, Cai, Lin and He introduced and investigated the concept of paratopological left Bol loops and proved some results of paratopological groups can be extended to paratopological left Bol loops. Banakh and Repov\v s \cite{Banakh} studied many generalized metric properties in rectifiable spaces and topological lops and showed that a rectifiable space $X$ is metrizable if and only if it is sequential, has
countable $cs^{*}$-character, and contains no closed copy of the Fr\'echet-Urysohn fan $S_{\omega}$. In \cite{Shen2020}, Shen introduced paratopological left-loops and showed that every weakly first-countable paratopological left-loop is first-countable. As we all known, a gyrogroup as an important type of non-associative algebra has many applications in Geometry and Physics, in particular, in the study of the $c$-ball of relativistically admissible velocities with the Einstein velocity addition, see \cite{UA1988,UA2002,UA2005,UA}. Therefore, topological gyrogroups are very important topological spaces which were posed by Atiponrat \cite{AW}. Clearly, every topological group is a topological gyrogroup and each topological gyrogroup is a rectifiable space. The readers may consult \cite{AW1,AW2020,BL,BL2,BL3,BX2022,BZX,BZX2,WAS2020,ZBX} for more details about topological gyrogroups.

In this paper, we show that each compact subset of a topological gyrogroup with an $\omega^{\omega}$-base is metrizable, which deduces that if $G$ is a topological gyrogroup with an $\omega^{\omega}$-base and is a $k$-space, then it is sequential. Moreover, we mainly research some weakly first-countable properties in feathered strongly topological gyrogroups. Indeed, for further study on M\"{o}bius gyrogroups, Bao and Lin posed the concept of strongly topological gyrogroups and showed that a strongly topological gyrogroup $G$ is feathered if and only if it contains a compact $L$-subgyrogroup $H$ such that the quotient space $G/H$ is metrizable. Based on the characterization of feathered strongly topological gyrogroups, they proved that each feathered strongly topological gyrogroup is paracompact. Also based on the characterization, we give some equivalent relationships of metrizability for strongly topological gyrogroups, such as every feathered strongly topological gyrogroup with countable $cs^{*}$-character is metrizable, and so on. It is also shown that for a feathered strongly topological gyrogroup $G$, $G$ has a compact resolution swallowing the compact sets if and only if $G$ contains a compact $L$-subgyrogroup $H$ such that the quotient space $G/H$ is a Polish space. Some problems about topological gyrogroups with countable $cs^{*}$-character are posed.

\section{Preliminary}

Throughout this paper, all topological spaces are assumed to be Hausdorff, unless otherwise is explicitly stated. Let $\mathbb{N}$ be the set of all positive integers and $\omega$ the first infinite ordinal. The readers may consult \cite{AA, E, linbook1, UA} for notation and terminology not explicitly given here. Next we recall some definitions and facts.

\begin{definition}\cite{AW}
Let $G$ be a nonempty set, and let $\oplus: G\times G\rightarrow G$ be a binary operation on $G$. Then the pair $(G, \oplus)$ is called a {\it groupoid}. A function $f$ from a groupoid $(G_{1}, \oplus_{1})$ to a groupoid $(G_{2}, \oplus_{2})$ is called a {\it groupoid homomorphism} if $f(x\oplus_{1}y)=f(x)\oplus_{2} f(y)$ for any elements $x, y\in G_{1}$. Furthermore, a bijective groupoid homomorphism from a groupoid $(G, \oplus)$ to itself will be called a {\it groupoid automorphism}. We write $\mbox{Aut}(G, \oplus)$ for the set of all automorphisms of a groupoid $(G, \oplus)$.
\end{definition}

\begin{definition}\cite{UA}
Let $(G, \oplus)$ be a groupoid. The system $(G,\oplus)$ is called a {\it gyrogroup}, if its binary operation satisfies the following conditions:

\smallskip
$(G1)$ There exists a unique identity element $0\in G$ such that $0\oplus a=a=a\oplus0$ for all $a\in G$.

\smallskip
$(G2)$ For each $x\in G$, there exists a unique inverse element $\ominus x\in G$ such that $\ominus x \oplus x=0=x\oplus (\ominus x)$.

\smallskip
$(G3)$ For all $x, y\in G$, there exists $\mbox{gyr}[x, y]\in \mbox{Aut}(G, \oplus)$ with the property that $x\oplus (y\oplus z)=(x\oplus y)\oplus \mbox{gyr}[x, y](z)$ for all $z\in G$.

\smallskip
$(G4)$ For any $x, y\in G$, $\mbox{gyr}[x\oplus y, y]=\mbox{gyr}[x, y]$.
\end{definition}

Notice that a group is a gyrogroup $(G,\oplus)$ such that $\mbox{gyr}[x,y]$ is the identity function for all $x, y\in G$. The definition of a subgyrogroup is given as follows.

\begin{definition}\cite{ST}
Let $(G,\oplus)$ be a gyrogroup. A nonempty subset $H$ of $G$ is called a {\it subgyrogroup}, denoted
by $H\leq G$, if $H$ forms a gyrogroup under the operation inherited from $G$ and the restriction of $\mbox{gyr}[a,b]$ to $H$ is an automorphism of $H$ for all $a,b\in H$.

\smallskip
Furthermore, a subgyrogroup $H$ of $G$ is said to be an {\it $L$-subgyrogroup}, denoted
by $H\leq_{L} G$, if $\mbox{gyr}[a, h](H)=H$ for all $a\in G$ and $h\in H$.
\end{definition}

\begin{definition}\cite{AW}
A triple $(G, \tau, \oplus)$ is called a {\it topological gyrogroup} if the following statements hold:

\smallskip
(1) $(G, \tau)$ is a topological space.

\smallskip
(2) $(G, \oplus)$ is a gyrogroup.

\smallskip
(3) The binary operation $\oplus: G\times G\rightarrow G$ is jointly continuous while $G\times G$ is endowed with the product topology, and the operation of taking the inverse $\ominus (\cdot): G\rightarrow G$, i.e. $x\rightarrow \ominus x$, is also continuous.
\end{definition}

Obviously, every topological group is a topological gyrogroup. However, every topological gyrogroup whose gyrations are not identically equal to the identity is not a topological group.

\begin{example}\cite[Example 3]{AW}
The Einstein gyrogroup with the standard topology is a topological gyrogroup but not a topological group.
\end{example}

Then we recall some weakly first-countable concepts which are important in the following researches.

\begin{definition}\cite{BT,GK,LPT}
A point $x$ of a topological space $X$ is said to have a {\it neighborhood $\omega^{\omega}$-base} or a {\it local $\mathfrak{G}$-base} if there exists a base of neighborhoods at $x$ of the form $\{U_{\alpha}(x):\alpha \in \mathbb{N}^{\mathbb{N}}\}$ such that $U_{\beta}(x)\subset U_{\alpha}(x)$ for all elements $\alpha \leq \beta$ in $\mathbb{N}^{\mathbb{N}}$, where $\mathbb{N}^{\mathbb{N}}$ consisting of all functions from $\mathbb{N}$ to $\mathbb{N}$ is endowed with the natural partial order, ie., $f\leq g$ if and only if $f(n)\leq g(n)$ for all $n\in \mathbb{N}$. The space $X$ is said to have an {\it $\omega^{\omega}$-base} or a {\it $\mathfrak{G}$-base} if it has a neighborhood $\omega^{\omega}$-base or a local $\mathfrak{G}$-base at every point $x\in X$.
\end{definition}

\begin{definition}
Let $X$ be a topological space.

$(1)$\, $X$ is called a {\it sequential space} \cite{FS} if for each non-closed subset $A\subseteq X$, there are a point $x\in X\setminus A$ and a sequence in $A$ converging to $x$ in $X$.

$(2)$\, $X$ is called a {\it Fr\'{e}chet-Urysohn space} \cite{FS} if for any subset $A\subseteq X$ and $x\in \overline{A}$, there is a sequence in $A$ converging to $x$ in $X$.

$(3)$\, $X$ is called an {\it $\alpha_{7}$-space} \cite{BZ}, if for every point $x\in X$ and each sheaf $\{S_{n}:n\in\omega\}$ with the vertex $x$, there exists a sequence converging to some point $y\in X$ which meets infinitely many sequences $S_{n}$.
\end{definition}

A family $\mathcal{N}$ of subsets of a topological space $X$ is called a {\it $cs^{*}$-network at a point $x\in X$} \cite{GMZ} if for each sequence $(x_{n})_{n\in \mathbb{N}}$ in $X$ convergent to $x$ and for each neighborhood $O_{x}$ of $x$ there is a set $N\in \mathcal{N}$ such that $x\in N\subset O_{x}$ and the set $\{n\in \mathbb{N}:x_{n}\in N\}$ is infinite.

Then we give the concept of $cs^{*}$-character of a topological gyrogroup.

\begin{definition}\cite[Theorem 3.7]{BZX2}
Let $G$ be a topological gyrogroup, the $cs^{*}$-character of $G$ is the least cardinality of $cs^{*}$-network at the identity element $0$ of $G$.
\end{definition}

\vspace{-0.2cm}

\section{Weakly first-countable properties of topological gyrogroups}

In this section, it is shown that each compact subset of a topological gyrogroup with an $\omega^{\omega}$-base is metrizable, which deduces that if $G$ is a topological gyrogroup with an $\omega^{\omega}$-base and is a $k$-space, then it is sequential. Moreover, for a feathered strongly topological gyrogroup $G$, based on the characterization of feathered strongly topological gyrogroups, that is, a strongly topological gyrogroup is feathered if and only if it contains a compact $L$-subgyrogroup $H$ such that the quotient space $G/H$ is metrizable, we give some equivalent relationships of metrizability for strongly topological gyrogroups, such as every feathered strongly topological gyrogroup with countable $cs^{*}$-character is metrizable.

\begin{definition}\cite[Definition 3.1]{GKL}
A family $\{K_{\alpha}:\alpha \in \mathbb{N}^{\mathbb{N}}\}$ of compact sets of a topological space $X$ is called a {\it compact resolution} if $X=\bigcup \{K_{\alpha}:\alpha \in \mathbb{N}^{\mathbb{N}}\}$ and $K_{\alpha}\subseteq K_{\beta}$ for all $\alpha \leq \beta$. In additionally, every compact set in $X$ is contained in some $K_{\alpha}$, we say that $\{K_{\alpha}:\alpha \in \mathbb{N}^{\mathbb{N}}\}$ {\it swallows the compact sets} of $X$.
\end{definition}

It is well-known that each Polish space $X$ has a compact resolution swallowing the compact sets of $X$. Moreover, it was proved in \cite[Theorem 3.3]{CRJP} that if $X$ is a metrizable topological space, then $X$ is a Polish space if and only if $X$ has a compact resolution swallowing the compact sets of $X$. Then, it follows from \cite[Proposition 3.3]{TVV} that each hemicompact topological space has a compact resolution swallowing the compact sets and the property of having a compact resolution swallowing the compact sets is closed-hereditary and is closed under countable products.

\begin{theorem}\label{3compact}
Let $G$ be a topological gyrogroup which has an $\omega^{\omega}$-base, $K$ an arbitrary compact subset of $G$. Then $K$ is metrizable.
\end{theorem}

\begin{proof}
By the hypothesis, let $\{U_{\alpha}:\alpha \in \mathbb{N}^{\mathbb{N}}\}$ be an open $\omega^{\omega}$-base in $G$. Without loss of generality, we assume that all sets $U_{\alpha}$ are symmetric. By \cite[Theorem 1]{CBOJ}, a compact space $K$ is metrizable if and only if $(K\times K)\setminus \Delta$ has a compact resolution swallowing its compact sets, where $\Delta =\{(x,x):x\in K\}$. Therefore, it suffices to show that the set $W=(K\times K)\setminus \Delta$ has a compact resolution which swallows its compact sets.

For each $\alpha\in \mathbb{N}^{\mathbb{N}}$, set $W_{\alpha}=\{(x,y)\in W,x\oplus (\ominus y)\not\in U_{\alpha}\}$. Then $W_{\alpha}$ is closed in $K\times K$, and hence it is compact for each $\alpha\in \mathbb{N}^{\mathbb{N}}$.

For each compact subset $C$ of $W$, $q(C)=\{x\oplus (\ominus y):(x,y)\in C\}$ is compact and does not contain the identity element $0$ of $G$. Since $\{U_{\alpha}:\alpha\in \mathbb{N}^{\mathbb{N}}\}$ is a local base at $0$, for some $\alpha\in \mathbb{N}^{\mathbb{N}}$, we obtain $U_{\alpha}\cap q(C)=\emptyset$. Then $C\subseteq W_{\alpha}$. Thus, $\{W_{\alpha}:\alpha\in \mathbb{N}^{\mathbb{N}}\}$ is a compact resolution swallowing the compact sets in $W$. We conclude that $K$ is metrizable.
\end{proof}

In \cite[Theorem 1.1]{Banakh}, Banakh showed that each non-metrizable sequential rectifiable space $X$ of countable $cs^{*}$-character contains a clopen rectifiable submetrizable $k_{\omega}$-subspace. Indeed, during the process of proof, it is not difficult to see that if $X$ is a non-metrizable sequential topological gyrogroup which has countable $cs^{*}$-character, then it contains an open and closed subgyrogroup which is a submetrizable $k_{\omega}$-space. For a topological space $X$, Chasco, Mart\'{i}n and Tarieladze in \cite[Lemma 1.5]{CMT} showed that if $X$ is sequential, then it is a $k$-space and if $X$ is a Hausdorff $k$-space and its compact subsets are sequential (in particular first countable or metrizable), then $X$ is sequential. Furthermore, it was proved in \cite[Theorem 3.8]{BZX2} that if a topological gyrogroup $G$ has an $\omega^{\omega}$-base, then it has countable $cs^{*}$-character. Therefore, by Theorem \ref{3compact} and these results, we obtain:

\begin{corollary}\label{k-sequential}
If a topological gyrogroup $G$ has an $\omega^{\omega}$-base, then the following conditions are equivalent.
\begin{enumerate}
\smallskip
\item $G$ is a $k$-space;

\smallskip
\item $G$ is sequential;

\smallskip
\item $G$ is metrizable or contains an open submetrizable $k_{\omega}$-subgyrogroup.
\end{enumerate}
\end{corollary}

In \cite[Theorem 3.5]{ZBX}, the authors showed that if $G$ is a sequential topological gyrogroup with an $\omega^{\omega}$-base, then $G$ has the strong Pytkeev property. Therefore, Corollary \ref{k-sequential} poses the following result directly.

\begin{theorem}\label{k-Pytkeev}
Let $G$ be a topological gyrogroup with an $\omega^{\omega}$-base. If $G$ is a $k$-space, then $G$ has the strong Pytkeev property.
\end{theorem}

Since for a Baire topological gyrogroup $G$, $G$ is metrizable if and only if it has the strong Pytkeev property, see \cite[Theorem 3.10]{ZBX}, Theorem \ref{k-Pytkeev} provides the following corollary.

\begin{corollary}
Let $G$ be a Baire topological gyrogroup. Then $G$ is metrizable if and only if $G$ is a $k$-space and has an $\omega^{\omega}$-base.
\end{corollary}

Theorem \ref{k-Pytkeev} also provides the other type of proof about the following result.

\begin{corollary}\cite[Corollary 3.6]{BZX2}
A topological gyrogroup $G$ is metrizable if and only if $G$ has an $\omega^{\omega}$-base and $G$ is also Fr\'echet-Urysohn.
\end{corollary}

\begin{proof}
The necessity is trivial, it suffices to prove the sufficiency.

Suppose that $G$ is a Fr\'echet-Urysohn topological gyrogroup and has an $\omega^{\omega}$-base. By Theorem \ref{k-Pytkeev}, $G$ has the strong Pytkeev property, hence has countable $cs^{*}$-character. It follows from \cite[Corollary 3.6]{BX2022} that every Fr\'echet-Urysohn topological gyrogroups with countable $cs^{*}$-character is metrizable.
\end{proof}

\begin{remark} A topological gyrogroup $G$ with an $\omega^{\omega}$-base has countable $cs^{*}$-character, see \cite[Theorem 3.8]{BZX2}, hence it is natural to consider the following question. If a topological gyrogroup $G$ is of countable $cs^{*}$-character and it is a $k$-space, then is $G$ sequential? Indeed, Shen in \cite[Example 4.5]{Shen2014} showed that there is a non-metrizable $snf$-countable topological group $X$ which is a $k$-space. Clearly, $X$ is not sequential. Furthermore, this example gives a negative answer to the question whether the $k$-property and sequentiality are equivalent for topological groups with countable $cs^{*}$-character posed in \cite{GKL} under Corollary 3.13.
\end{remark}

Then, we introduce the concept of strongly topological gyrogroups, which was first posed by Bao and Lin in \cite{BL}, and we investigate some weakly first-countable properties in feathered strongly topological gyrogroups.

\begin{definition}{\rm (\cite{BL})}\label{d11}
Let $G$ be a topological gyrogroup. We say that $G$ is a {\it strongly topological gyrogroup} if there exists a neighborhood base $\mathscr U$ of $0$ such that, for every $U\in \mathscr U$, $\mbox{gyr}[x, y](U)=U$ for any $x, y\in G$. For convenience, we say that $G$ is a strongly topological gyrogroup with neighborhood base $\mathscr U$ of $0$.
\end{definition}

For each $U\in \mathscr U$, we can set $V=U\cup (\ominus U)$. Then, $$\mbox{gyr}[x,y](V)=\mbox{gyr}[x, y](U\cup (\ominus U))=\mbox{gyr}[x, y](U)\cup (\ominus \mbox{gyr}[x, y](U))=U\cup (\ominus U)=V,$$ for all $x, y\in G$. Obviously, the family $\{U\cup(\ominus U): U\in \mathscr U\}$ is also a neighborhood base of $0$. Therefore, we may assume that $U$ is symmetric for each $U\in\mathscr U$ in Definition~\ref{d11}. Moreover, in the classical M\"{o}bius, Einstein, or Proper Velocity gyrogroups, we know that gyrations are indeed special rotations, however for an arbitrary gyrogroup, gyrations belong to the automorphism group of $G$ and need not be necessarily rotations.

In \cite{BL}, the authors proved that there is a strongly topological gyrogroup which is not a topological group.

\begin{example}\cite[Example 3.1]{BL}
The M\"{o}bius gyrogroup with the standard topology is a strongly topological gyrogroup but not a topological group.
\end{example}

A topological gyrogroup $G$ is {\it feathered} if it contains a non-empty compact set $K$ of countable character in $G$. It was proved in \cite[3.1 E(b) and 3.3 H(a)]{E} that every locally compact topological gyrogroup is feathered. Moreover, by \cite[Theorem 3.14]{BL}, we know that a strongly topological gyrogroup $G$ is feathered if and only if it contains a compact $L$-subgyrogroup $H$ such that the quotient space $G/H$ is metrizable.

\begin{theorem}\label{feathered}
Let $G$ be a feathered strongly topological gyrogroup. Then $G$ has an $\omega^{\omega}$-base if and only if $G$ is metrizable.
\end{theorem}

\begin{proof}
Suppose that $G$ is a strongly topological gyrogroup and has an $\omega^{\omega}$-base. Then $G$ contains a compact $L$-subgyrogroup $H$ such that the quotient space $G/H$ is metrizable. It follows from Theorem \ref{3compact} that the subgyrogroup $H$ is metrizable. Since each compact subset of a Hausdorff space is closed, it is clear that $H$ is a closed $L$-subgyrogroup of $G$. Then, by \cite[Corollary 4.3]{BZX}, if $G$ is a topological gyrogroup and $H$ is a closed $L$-subgyrogroup of $G$ and if the spaces $H$ and $G/H$ are metrizable, then the space $G$ is also metrizable. Therefore, we obtain that $G$ is a metrizable space.
\end{proof}

\begin{corollary}
Let $G$ be a locally compact strongly topological gyrogroup. Then $G$ has an $\omega^{\omega}$-base if and only if $G$ is metrizable.
\end{corollary}

Since every topological gyrogroup with an $\omega^{\omega}$-base has countable $cs^{*}$-character, it is natural to pose the following question.

\begin{question}\label{question-cs}
Let $G$ be a feathered strongly topological gyrogroup with countable $cs^{*}$-character. Is $G$ metrizable?
\end{question}

Then, we give a affirmative answer to Question \ref{question-cs}. We note that Uspenski\v\i \cite{UVV, UVV89} proved that compact rectifiable spaces are dyadic. Since every topological gyrogroup is a rectifiable space, it is trivial that each compact topological gyrogroup is dyadic. Moreover, Banakh and Zdomskyy in \cite[Proposition 7]{BZ} claimed that a dyadic compactum is metrizable if and only if it has countable $cs^{*}$-character.

\begin{proposition}
Every compact topological gyrogroup with countable $cs^{*}$-character is metrizable.
\end{proposition}

\begin{theorem}\label{csf-feath}
Every feathered strongly topological gyrogroup with countable $cs^{*}$-character is metrizable.
\end{theorem}

\begin{proof}
Since $G$ is a feathered strongly topological gyrogroup, there exists a compact $L$-subgyrogroup $H$ of $G$ such that the quotient space $G/H$ is metrizable. By the hypothesis, $G$ has countable $cs^{*}$-character, then $H$ also has countable $cs^{*}$-character, which deduces that the compact subgyrogroup $H$ with countable $cs^{*}$-character is metrizable, and it follows from \cite[Corollary 4.3]{BZX} that $G$ is metrizable.
\end{proof}

\begin{corollary}
Every locally compact strongly topological gyrogroup with countable $cs^{*}$-character is metrizable.
\end{corollary}

\begin{theorem}
Let $G$ be a strongly topological gyrogroup. Then the following conditions are equivalent:
\begin{enumerate}
\smallskip
\item $G$ is metrizable;

\smallskip
\item $G$ is Fr\'echet-Urysohn and has countable $cs^{*}$-character;

\smallskip
\item $G$ is Fr\'echet-Urysohn and has an $\omega^{\omega}$-base;

\smallskip
\item $G$ is feathered and has countable $cs^{*}$-character;

\smallskip
\item $G$ is feathered and has an $\omega^{\omega}$-base.
\end{enumerate}
\end{theorem}

In the research of rectifiable spaces, Banakh and Repov\v s \cite[Lemma 5.1]{Banakh} showed that suppose that $G$ is a topological lop and $F\subseteq G$ is a subset containing the unit $e$ of $G$, then put $F_{1}=F$ and $F_{n+1}=F_{n}^{-1}F_{n}$ for $n\in \mathbb{N}$. If $F$ is a sequential space containing no closed topological copy of the Fr\'echet-Urysohn fan $S_{\omega}$ and each space $F_{n}$, $n\in \mathbb{N}$, has a countable $cs^{*}$-network at $e$, then $F$ has a countable $sb$-network at $e$; if $F$ is sequential and each space $F_{n},n\in \mathbb{N}$, has countable $sb$-network at $e$, then $F$ is first-countable at $e$. Then, Shen \cite[Proposition 2.6 and Theorem 2.7]{Shen2020} showed that every paratopological left-loop with $sb$-network is $sof$-countable and a sequential, regular paratopological left-loop $G$ with countable $cs^{*}$-network is first-countable if and only if $G$ contains no closed copy of $S_{\omega}$. These results in both of two articles can obtain the following results immediately.

\begin{proposition}\label{csf-snf}
If $G$ is a sequential topological gyrogroup with countable $cs^{*}$-network containing no closed copy of $S_{\omega}$, then $G$ has countable $sb$-network.
\end{proposition}

\begin{proposition}\label{snf-sof}
Every topological gyrogroup with countable $sb$-network is $sof$-countable.
\end{proposition}

\begin{theorem}
A strongly topological gyrogroup $G$ is metrizable if and only if $G$ is a $k$-space of countable pseudocharacter with countable $sb$-network.
\end{theorem}

\begin{proof}
The necessity is trivial, it suffices to claim the sufficiency.

Let a strongly topological gyrogroup $G$ be a $k$-space of countable pseudocharacter with countable $sb$-network. It follows from \cite[Theorem 4.3]{BL1} that every strongly topological gyrogroup with countable pseudocharacter is submetrizable. We obtain that every compact subset of $G$ is metrizable. Since $G$ is a $k$-space, it is easy to see that $G$ is sequential. Indeed, a space $X$ is first-countable if and only if $X$ is sequential and $sof$-countable, which deduces that $G$ is first-countable by Proposition \ref{snf-sof}, hence $G$ is metrizable.
\end{proof}

Recall that a continuous mapping $q:G\rightarrow H$ is called {\it compact-covering} if for every compact subset $K$ of $H$ there exists a compact subset $C$ of $G$ such that $q(C)=K$. Indeed, it was claimed in \cite[Theorem 3.8]{BL} that if $G$ is a topological gyrogroup and $H$ is a compact L-subgyrogroup of $G$, then the quotient mapping $\pi$ of $G$ onto the quotient space $G/H$ is perfect. Furthermore, if $f:X\rightarrow Y$ is a perfect mapping, then for every compact subspace $Z\subseteq Y$, the inverse image $f^{-1}(Z)$ is compact by \cite[Theorem 3.7.2]{E}. Therefore, if $H$ is a compact $L$-subgyrogroup of a topological gyrogroup $G$, then the quotient mapping $\pi$ of $G$ onto the quotient space $G/H$ is a compact covering mapping.

\begin{theorem}\label{Polish}
If $G$ is a feathered strongly topological gyrogroup, the followings are equivalent.

\begin{enumerate}
\smallskip
\item $G$ has a compact resolution swallowing the compact sets of $G$;

\smallskip
\item $G$ has a compact $L$-subgyrogroup $H$ such that the quotient space $G/H$ is a Polish space.
\end{enumerate}

\smallskip
If (1) holds, then $G$ is $\check{C}$ech-complete.
\end{theorem}

\begin{proof}
(1)$\Rightarrow$(2). Let $\pi$ be the natural homomorphism from $G$ to its quotient topology on $G/H$. Since $G$ is a feathered strongly topological gyrogroup, it follows from \cite[Lemma 3.14]{BL} that $G$ contains a compact $L$-subgyrogroup $H$ such that $G/H$ is metrizable. Let $G$ have a compact resolution swallowing the compact sets of $G$, say $\mathcal{K}=\{K_{\alpha}:\alpha\in \mathbb{N}^{\mathbb{N}}\}$. Then put $\mathcal{K}'=\{\pi (K_{\alpha}):\alpha\in \mathbb{N}^{\mathbb{N}}\}$. Then $\mathcal{K}'$ swallows the compact subsets of $G/H$. Indeed, if $K'$ is compact in $G/H$, then $\pi^{-1}(K')$ is compact in $G$. So there exists $\alpha\in \mathbb{N}^{\mathbb{N}}$ such that $\pi^{-1}(K')\subseteq K_{\alpha}$ and hence $K'\subseteq \pi (K_{\alpha})$. We know that $G/H$ is a Polish space by \cite[Theorem 3.3]{CRJP}.

(2)$\Rightarrow$(1). Since the space $G/H$ is a Polish space, it follows from \cite[Theorem 3.3]{CRJP} that $G/H$ has a compact resolution swallowing the compact sets of $G/H$, say $\mathcal{K}'=\{K'_{\alpha}:\alpha\in \mathbb{N}^{\mathbb{N}}\}$. For every $\alpha\in \mathbb{N}^{\mathbb{N}}$, put $K_{\alpha}=\pi^{-1}(K'_{\alpha})$. Then $K_{\alpha}$ is a compact subset of $G$. Hence, $\mathcal{K}=\{K_{\alpha}:\alpha\in \mathbb{N}^{\mathbb{N}}\}$ is a compact resolution. Let $C$ be a compact subset of $G$. Then there exists $\alpha\in \mathbb{N}^{\mathbb{N}}$ such that $\pi (C)\subseteq K'_{\alpha}$. Therefore, $C\subseteq K_{\alpha}$, and hence $\mathcal{K}$ swallows the compact sets of $G$. We conclude that $G$ has a compact resolution swallowing the compact sets of $G$.

By \cite[Theorem 3.17]{BL}, a strongly topological gyrogroup $G$ is $\check{C}$ech-complete if and only if it contains a compact $L$-subgyrogroup $H$ such that the quotient space $G/H$ is metrizable by a complete metric. Since each Polish space is a complete metric space, we know that $G$ is $\check{C}$ech-complete.
\end{proof}

\begin{proposition}
Let $G$ be a topological gyrogroup and have a compact resolution swallowing the compact sets of $G$. If $q:G\rightarrow H$ is a quotient compact-covering map, then $H$ also has a compact resolution swallowing the compact sets of $H$.
\end{proposition}

\begin{proof}
Indeed, the result is trivial. If $\{K_{\alpha}:\alpha\in \mathbb{N}^{\mathbb{N}}\}$ is a compact resolution swallowing the compact sets of $G$, it is clear that $\{q(K_{\alpha}):\alpha\in \mathbb{N}^{\mathbb{N}}\}$ is a compact resolution swallowing the compact sets of $H$.
\end{proof}

\begin{proposition}
Let $G$ be a topological gyrogroup and a $k$-space. If $G$ has an $\omega^{\omega}$-base and also has a compact resolution swallowing the compact sets of $G$, then $G$ is either a Polish space or contains a submetrizable $k_{\omega}$-subgyrogroup.
\end{proposition}

\begin{proof}
Since topological gyrogroup $G$ is a $k$-space and has an $\omega^{\omega}$-base, we have that $G$ is metrizable or contains an open submetrizable $k_{\omega}$-subgyrogroup. Moreover, by \cite[Theorem 3.3]{CRJP}, in a metrizable space $X$, $X$ is a Polish space if and only if $X$ has a compact resolution swallowing the compact sets of $X$. Therefore, if $G$ is metrizable, we know that $G$ is a Polish space.
\end{proof}

In Theorems \ref{feathered}, \ref{csf-feath} and \ref{Polish}, it is clear that the characterization of feathered strongly topological gyrogroup playes an important role in the proof. However, we do not know whether the characterization of feathered holds in topological gyrogroups. If it holds in topological gyrogroups, all of Theorems \ref{feathered}, \ref{csf-feath} and \ref{Polish} can be extended to topological gyrogroups immediately.

\begin{question}
If $G$ is a feathered topological gyrogroup, is there a compact $L$-subgyrogroup of $G$ such that the quotient space $G/H$ metrizable?
\end{question}

A space $X$ is called {\it hemicompact} if $X=\bigcup_{n\in \mathbb{N}}X_{n}$, where $X_{n}$ is compact for any $n\in \mathbb{N}$ and for any compact $K\subseteq X$, there is $n\in \mathbb{N}$ such that $K\subseteq X_{n}$.

\begin{corollary}
Let $G$ be a locally compact strongly topological gyrogroup. Then $G$ has a compact resolution swallowing the compact sets of $G$ if and only if $G$ is hemicompact space.
\end{corollary}

\begin{proof}
Since each hemicompact topological space has a compact resolution swallowing the compact sets, it suffices to claim the necessity.

Suppose that $G$ has a compact resolution swallowing the compact sets of $G$. Since each locally compact topological gyrogroup is feathered, by Theorem \ref{Polish}, $G$ contains a compact $L$-subgyrogroup $H$ such that the locally compact space $G/H$ is second countable. Therefore, $G/H$ is hemicompact, and $G$ is also hemicompact.
\end{proof}

\begin{proposition}
Every Fr\'echet-Urysohn hemicompact topological gyrogroup is locally compact.
\end{proposition}

\begin{proof}
Let $G=\bigcup_{n\in \mathbb{N}}K_{n}$, where $\{K_{n}\}_{n}$ is an increasing sequence of compact subsets of $K$ containing the identity element $0$ such that every compact set in $G$ is contained in some $K_{n}$. Then we can find $n\in \mathbb{N}$ such that $K_{n}$ is a neighborhood of $0$. Suppose on the contrary that there is no $n$ such that $K_{n}$ is a neighborhood of $0$ for each $n\in \mathbb{N}$. Then for each $n\in \mathbb{N}$ and each neighborhood $U$ of $0$, there exists $x_{U,n}\in U\setminus K_{n}$. For each $n\in \mathbb{N}$, set $B_{n}=\{x_{U,n}:U\mbox{ is an open neighborhood of 0}\}$. Then $0\in \overline{B_{n}}$. Since $G$ is Fr\'echet-Urysohn, for each $n\in \mathbb{N}$, we can find an open neighborhood sequence $\{U_{n}(k)\}_{k}$ of $0$ such that $x_{U_{n}(k),n}\rightarrow 0$ at $k\rightarrow \infty$. Since every Fr\'echet-Urysohn topological gyrogroup is a strong $\alpha_{4}$-space by \cite[Lemma 3.3]{BZX2}, there exists strictly increasing sequences $(k_{p})_{p}$ and $(n_{p})_{p}$ such that $x_{U_{n_{p}}(k_{p}),n_{p}}\rightarrow 0$ at $p\rightarrow \infty$. Since the set $B=\{x_{U_{n_{p}}(k_{p}),n_{p}}:p\in \mathbb{N}\}\cup \{0\}$ is compact in $G$, we can find $m\in \mathbb{N}$ such that $B\subseteq K_{m}$, which is a contradiction . Therefore, $G$ is locally compact.
\end{proof}

\section{metrizability of strongly topological gyrogroups}

In \cite{AW}, Atiponrat posed a question that is it true that the first-countability axiom implies that $G$ is metrizable for a topological gyrogroup $G$? Then Cai, Lin and He showed that every topological gyrogroup is a rectifiable space, which deduces that every first-countable (strongly) topological gyrogroup is metrizable. Indeed, Alexandra S. Gul'ko \cite[Theorem 3.2]{Alexan} proved that every first-countable $T_{0}$ rectifiable space is metrizable by the tool of strong development. Here, we give a direct construction to show that every first-countable strongly topological gyrogroup is metrizable.

\begin{lemma}\label{3.3.10}\cite[Lemma 3.12]{BL}
Let $G$ be a strongly topological gyrogroup with the symmetric neighborhood base $\mathscr{U}$ at $0$, and let $\{U_{n}: n\in\mathbb{N}\}$ and $\{V(m/2^{n}): n, m\in\mathbb{N}\}$ be two sequences of open neighborhoods satisfying the following conditions (1)-(5):

\smallskip
(1) $U_{n}\in\mathscr{U}$ for each $n\in \mathbb{N}$.

\smallskip
(2) $U_{n+1}\oplus U_{n+1}\subseteq U_{n}$, for each $n\in \mathbb{N}$.

\smallskip
(3) $V(1)=U_{0}$;

\smallskip
(4) For any $n\geq 1$, put $$V(1/2^{n})=U_{n}, V(2m/2^{n})=V(m/2^{n-1})$$ for $m=1,...,2^{n-1}$, and $$V((2m+1)/2^{n})=U_{n}\oplus V(m/2^{n-1})=V(1/2^{n})\oplus V(m/2^{n-1})$$ for each $m=1,...,2^{n-1}-1$;

\smallskip
(5) $V(m/2^{n})=G$ when $m>2^{n}$;

\smallskip
Then there exists a prenorm $N$ on $G$ that satisfies the following conditions:

\smallskip
(a) for any fixed $x, y\in G$, we have $N(\mbox{gyr}[x,y](z))=N(z)$ for any $z\in G$;

\smallskip
(b) for any $n\in \mathbb{N}$, $$\{x\in G: N(x)<1/2^{n}\}\subseteq U_{n}\subseteq\{x\in G: N(x)\leq 2/2^{n}\}.$$
\end{lemma}

\begin{theorem}\label{new-metric}
Every first-countable strongly topological gyrogroup is metrizable.
\end{theorem}

\begin{proof}
Let $G$ be a strongly topological gyrogroup with a symmetric neighborhood base $\mathscr U$. Since $G$ is first-countable, put $\{W_{n}:n\in \mathbb{N}\}$ a countable base at the identity element $0$. By induction, we obtain a sequence $\{U_{n}:n\in \mathbb{N}\}\subseteq \mathscr{U}$ such that $U_{n}\subseteq W_{n}$ and $U_{n+1}\oplus U_{n+1}\subseteq U_{n}$, for each $n\in \mathbb{N}$. It is easy to see that $\{U_{n}:n\in \mathbb{N}\}$ is also a base of $G$ at $0$. By Lemma \ref{3.3.10}, there exists a continuous prenorm $N$ on $G$ which satisfies $$N(\mbox{gyr}[x, y](z))=N(z)$$ for any $x, y, z\in G$ and $$\{x\in G: N(x)<1/2^{n}\}\subseteq U_{n}\subseteq \{x\in G: N(x)\leq 2/2^{n}\},$$ for each integer $n\geq 0$. Put $B_{N}(\varepsilon)=\{x\in G:N(x)<\varepsilon\}$, where $\varepsilon$ is a positive number. It is easy to see that $B_{N}(1/2^{n})$ also forms a base of $G$ at $0$.

Now, for arbitrary $x$ and $y$ in $G$, put $\varrho _{N}(x, y)=N(\ominus x\oplus y)$. Let us show that $\varrho _{N}$ is a metric on $G$.

\smallskip
(1) Clearly, $\varrho _{N}(x, y)=N(\ominus x\oplus y)\geq 0$, for every $x, y\in G$. At the same time, $\varrho _{N}(x, x)=N(0)=0$, for each $x\in G$. Assume that $$\varrho _{N}(x, y)=N(\ominus x\oplus y)=0.$$ Then, for each $n\in\mathbb{N}$, $$\ominus x\oplus y\in \{x\in G: N(x)<1/2^{n}\}\subseteq U_{n}.$$ Since $\{0\}=\bigcap _{n\in \mathbb{N}}U_{n}$, it follows that $\ominus x\oplus y=0$, that is, $x=y$.

\smallskip
(2) For every $x, y\in G$, $\varrho _{N}(y, x)=N(\ominus y\oplus x)=N(gyr[\ominus y,x](\ominus x\oplus y))=N(\ominus x\oplus y)=\varrho _{N}(x, y)$.

\smallskip
(3) For every $x, y, z\in G$, it follows from \cite[Theorem 2.11]{UA2005} that
\begin{eqnarray}
\varrho _{N}(x, y)&=&N(\ominus x\oplus y)\nonumber\\
&=&N((\ominus x\oplus z)\oplus \mbox{gyr}[\ominus x, z](\ominus z\oplus y))\nonumber\\
&\leq&N(\ominus x\oplus z)+N(\mbox{gyr}[\ominus x, z](\ominus z\oplus y))\nonumber\\
&=&N(\ominus x\oplus z)+N(\ominus z\oplus y)\nonumber\\
&=&\varrho _{N}(x, z)+\varrho _{N}(z, y)\nonumber
\end{eqnarray}

Thus, $\varrho _{N}$ is a metric on $G$.

Since $B_{N}(1/2^{n})$ forms a base of $G$ at $0$ and $G$ is homogeneous, for each $x\in G$, $B_{N}(1/2^{n})\oplus x$ constitutes a base of $G$ at $x$. Therefore, it is easy to see that the topology generated by metric $\varrho _{N}$ is coincide with the original topology of $G$. Hence, $G$ is metrizable.
\end{proof}

However, we do not know whether each topological gyrogroup has the similar result like Lemma \ref{3.3.10}, therefore, we can not give the direct construction that every first-countable topological gyrogroup is metrizable. Moreover, we pose the following questions.

\begin{question}
Let $G$ be a metrizable (strongly) topological gyrogroup and $H$ a closed $L$-subgyrogroup of $G$. Is the quotient space $G/H$ metrizable?
\end{question}

\begin{remark}
It was posed a question in \cite{BL} that if $G$ is a (strongly) topological gyrogroup with a countable pseudocharacter, is $G$ submetrizable? then Bao and Lin gave an affirmative answer to this question in \cite[Theorem 4.3]{BL1} by constructing a metric $\varrho _{N}(x, y)=N(\ominus x\oplus y)+N(\ominus y\oplus x)$ when $G$ is a strongly topological gyrogroup. Notice that the proof of Theorem \ref{new-metric} can be applied to show that every strongly topological gyrogroup with a countable pseudocharacter is submetrizable, that is, the metric in \cite[Theorem 4.3]{BL1} can be replaced by $\varrho _{N}(x, y)=N(\ominus x\oplus y)$.
\end{remark}

A subgyrogroup $H$ of a topological gyrogroup $G$ is called {\it inner (outer) neutral} if for every open neighborhood $U$ of $0$ in $G$, there exists an open neighborhood $V$ of $0$ such that $H\oplus V\subseteq U\oplus H$ ($V\oplus H\subseteq H\oplus U$).

\begin{question}
Let $G$ be a (strongly) topological gyrogroup and $H$ a closed inner neutral $L$-subgyrogroup of $G$. If the quotient space $G/H$ is first-countable, is it metrizable?
\end{question}

\begin{question}
Let $G$ be a feathered (strongly) topological gyrogroup and $H$ a closed $L$-subgyrogroup of $G$. If the quotient space $G/H$ is first-countable, is it metrizable?
\end{question}

\end{document}